\documentclass[11pt, a4paper]{article}
\usepackage{times}
\usepackage{a4wide}
\usepackage[dvips, hyperindex]{hyperref}
\usepackage[british]{babel}
\usepackage{enumerate, longtable}
\usepackage{amsmath, amscd, amsfonts, amssymb, latexsym, theorem, comment, graphicx}
\usepackage{xypic}
\usepackage[T1]{fontenc}
\usepackage[latin1]{inputenc}

{\theorembodyfont{\itshape}   \newtheorem{thm}{Theorem}[section]}
{\theorembodyfont{\itshape}   }
{\theorembodyfont{\itshape}   \newtheorem{defi}[thm]{Definition}}
{\theorembodyfont{\rmfamily}   \newtheorem{rem}[thm]{Remark}}
{\theorembodyfont{\itshape}   }
{\theorembodyfont{\itshape}   \newtheorem{prop}[thm]{Proposition}}
{\theorembodyfont{\itshape}   \newtheorem{cor}[thm]{Corollary}}
{\theorembodyfont{\itshape}   }
{\theorembodyfont{\itshape}   }
{\theorembodyfont{\itshape}   }
{\theorembodyfont{\itshape}   \newtheorem{hyp}[thm]{Hypothesis}}
{\theorembodyfont{\itshape}   \newtheorem{problem}[thm]{Problem}}
{\theorembodyfont{\itshape}   \newtheorem{question}[thm]{Question}}
{\theorembodyfont{\itshape}   }

\newcommand{\End}{{\rm End}}
\newcommand{\Hom}{{\rm Hom}}

\newcommand{\GL}{\mathrm{GL}}
\newcommand{\SL}{\mathrm{SL}}

\newcommand{\old}{\mathrm{old}}

\DeclareMathOperator{\Gal}{Gal}

\DeclareMathOperator{\Frob}{Frob}
\DeclareMathOperator{\Tr}{Tr}

\DeclareMathOperator{\Det}{det}

\newcommand{\Frac}{\mathrm{Frac}}

\newcommand{\pf}{{\bf Proof. }}

\newcommand{\qed}{\hspace* {.5cm} \hfill $\Box$}

\newcommand{\cO}{\mathcal{O}}

\newcommand{\fm}{\mathfrak{m}}

\newcommand{\fp}{\mathfrak{p}}

\newcommand{\fP}{\mathfrak{P}}

\newcommand{\CC}{\mathbb{C}}

\newcommand{\FF}{\mathbb{F}}

\newcommand{\NN}{\mathbb{N}}

\newcommand{\QQ}{\mathbb{Q}}

\newcommand{\TT}{\mathbb{T}}

\newcommand{\ZZ}{\mathbb{Z}}

\newcommand{\Qbar}{\overline{\QQ}}
\newcommand{\Fbar}{\overline{\FF}}

\newcommand{\Zbar}{{\overline{\ZZ}}}

\newcommand{\rhobar}{\overline{\rho}}

\newcommand{\ilim}[1]{\underset{#1}{\underrightarrow{\lim} \;}}

\begin{document}

\selectlanguage{british}

\title{Computing Congruences of\\
Modular Forms and Galois Representations\\ Modulo Prime Powers}
\author{Xavier Taixés i Ventosa\footnote{Universitat Pompeu Fabra,
Departament d'Economia i Empresa, Ramon Trias Fargas 25-27, 08005
Barcelona\newline xavier.taixes@upf.edu} \ and Gabor
Wiese\footnote{Universität Duisburg-Essen, Institut für
Experimentelle Mathematik, Ellernstraße 29, 45326 Essen,
Germany\newline gabor@pratum.net, http://maths.pratum.net/ }}
\maketitle

\begin{abstract}
This article starts a computational study of congruences of modular forms
and modular Galois representations modulo prime powers.
Algorithms are described that compute the maximum integer modulo which two
monic coprime integral polynomials have a root in common in a sense
that is defined.
These techniques are applied to the study of congruences of modular forms
and modular Galois representations modulo prime powers.
Finally, some computational results with implications
on the (non-)liftability of modular forms modulo prime powers
and possible generalisations of level raising are presented.

\end{abstract}

\section{Introduction}

Congruences of modular forms modulo a prime~$\ell$ and -- from a different point of view --
modular forms over~$\Fbar_\ell$ play an important role in modern Arithmetic Geometry. The most
prominent recent example is Serre's modularity conjecture, which has just become a theorem
of Khare, Wintenberger and Kisin. We particularly mention the various
techniques for {\em Level Raising} and {\em Level Lowering} modulo~$\ell$ that were already
crucial for Wiles's proof of Fermat's Last Theorem.

Motivated by this, it is natural to study congruences modulo~$\ell^n$ of
modular forms and Galois representations. However, as working
over non-factorial and non-reduced rings like $\ZZ/\ell^n\ZZ$ introduces many extra difficulties,
one is led to first approach this subject from an algorithmic and computational point of view,
which is the topic of this article.

We introduce a definition of when two algebraic integers $a$, $b$ are congruent
modulo~$\ell^n$. Our definition, which might appear non-standard at first,
was forced upon us by three requirements: Firstly, we want it to be independent of
any choice of number field containing $a$,~$b$.
Secondly, in the special case $n=1$ a congruence modulo~$\ell$
should come down to an equality in a finite field. Finally, if $a$, $b$
lie in some number field~$K$ that is unramified at~$\ell$, then a congruence of $a$
and $b$ modulo~$\ell^n$ should be a congruence modulo~$\lambda^n$, where $\lambda$ is a prime
dividing~$\ell$ in~$K$.

Since algebraic integers are -- up to Galois conjugacy -- most conveniently represented
by their minimal polynomials, we address the problem of determining for which
prime powers~$\ell^n$ two coprime monic integral polynomials have zeros
which are congruent modulo~$\ell^n$.
We prove that a certain number, called the reduced discriminant or -- in our language --
the congruence number of the two polynomials, in all cases gives a good upper bound
and in favourable cases completely solves this problem.
In the cases when the congruence number is insufficient, we use a method based
on the Newton polygon of the polynomial whose roots are the differences of the
roots of the polynomials we started with.

With these tools at our disposal, we target the problem of computing congruences
modulo~$\ell^n$ between two Hecke eigenforms. Since our motivation comes from arithmetic,
especially from Galois representations, our main interest is in Hecke eigenforms. It quickly
turns out, however, that there are several possible well justified notions of
Hecke eigenforms modulo~$\ell^n$. We present two, which we call {\em strong} and {\em weak}.
The former can be thought of as reductions modulo~$\ell^n$ of $q$-expansions of holomorphic
normalised Hecke eigenforms; the latter can be understood as linear combinations of
holomorphic modular forms, which are in general not eigenforms, but whose reduction modulo~$\ell^n$
becomes an eigenform (our definition is formulated in a different way, but can be interpreted
to mean this).
We observe that Galois representations to $\GL_2(R)$,
where $R$ is an extension of~$\ZZ/\ell^n\ZZ$ in the
sense of Section~\ref{sec-cong}, can be attached to both weak and strong Hecke eigenforms
(under the condition of residual absolute irreducibility).

Modular forms can be represented by their $q$-expansions (e.g.\ in $\ZZ/\ell^n\ZZ$), i.e.\
by power series. For computational purposes, such as uniquely identifying a modular
form and comparing two modular forms, it is essential that already a finite segment of a
certain length of the $q$-expansions suffices. We notice that a sufficient length
is provided by the so-called Sturm bound,
which is the same modulo~$\ell^n$ as in characteristic~$0$.

The computational problem that we are mostly interested in is to determine congruences
modulo~$\ell^n$ between two newforms, i.e.\ equalities between strong Hecke eigenforms modulo~$\ell^n$.
This problem is perfectly suited for applying our methods of determining congruences
modulo~$\ell^n$ of zeros of integral polynomials.
The reason for this is that the Fourier coefficient~$a_p$ of
a normalised Hecke eigenform is a zero of the characteristic polynomial of the Hecke operator
$T_p$ acting on a suitable integral modular symbols space (see e.g.\ \cite{S} or \cite{HeckeMS}).
Thus, in order to determine the prime powers modulo which two newforms are congruent, we compute
the congruences between the roots of these characteristic polynomials for a suitable
number of~$p$.
One important point deserves to be mentioned here: If the two newforms that we want to compare
do not have the same levels (but the same weights),
one cannot expect that they are congruent at all primes; a different
behaviour is to be expected at primes dividing the levels. We address this problem by
applying the usual degeneracy maps `modulo~$\ell^n$' in order to land in the same level.
All these considerations lead to an algorithm, which we sketch.
We point out that this algorithm is much faster than the (naive) one which works with the
coefficients of the modular forms as algebraic integers in a (necessarily big)
number field.

We implemented the algorithm and performed many computations which led to observations
that we consider very interesting.
Some of the results are reported upon in Section~\ref{sec-examples}. We are planning
to investigate questions like `Level Raising' in more detail in a subsequent work.
We remark that the algorithm was already used in~\cite{DT} to determine
some numerical examples satisfying the main theorem of that article.

\subsection*{Acknowledgements}

X.T.\ would like to thank Gerhard Frey for suggesting the subject of the article as
PhD project.
G.W.\ would like to thank Frazer Jarvis, Lara Thomas, Christophe Ritzenthaler,
Ian Kiming and, in particular, Gebhard Böckle for enlightening discussions
and e-mail exchanges relating to
the subject of this article, as well as Kristin Lauter for pointing out the
article~\cite{Pohst}. Special thanks are due to Michael Stoll for suggesting
the basic idea of one algorithm, as well as to one of the referee for also suggesting it
together with many other improvements in notation and presentation.
Thanks are also due to the second referee for pointing out that there should be
a relation to the paper~\cite{ARS}.

Both authors acknowledge partial support by the European Research Training Network
{\em Galois Theory and Explicit Methods} MRTN-CT-2006-035495.
G.~W. also acknowledges partial support by the Sonderforschungsbereich Transregio 45
of the Deutsche Forschungsgemeinschaft.

\subsection*{Notation}

We introduce some standard notation to be used throughout.
In the article $\ell$ and $p$ always refer to prime numbers.
By an $\ell$-adic field we shall understand a finite field extension of~$\QQ_\ell$.
We fix algebraic closures $\Qbar$ of~$\QQ$ and $\Qbar_\ell$ of~$\QQ_\ell$.
By $\Zbar$ and $\Zbar_\ell$ we denote the integers of $\Qbar$ and $\Qbar_\ell$, respectively.
If $K$ is either a number field or a local field, then $\cO_K$ denotes
its ring of integers. In the latter case, $\pi_K$ denotes a uniformiser, i.e.\
a generator of the maximal ideal of~$\cO_K$, and $v_K$ is the
valuation satisfying $v_K(\pi_K)=1$. Moreover, $v_\ell$ denotes the
valuation on $K$ and on $\Qbar_\ell$ normalised such that $v_\ell(\ell)=1$.

\section{Congruences modulo $\ell^n$}\label{sec-cong}

In this section we give our definition of {\em congruences modulo~$\ell^n$} for algebraic
and $\ell$-adic integers and discuss how to compute them.

\subsection{Definition}\label{subsec:def}

Since a question on congruences is a local question, we place ourselves
in the set-up of $\ell$-adic fields.
Let $\alpha,\beta \in \Zbar_\ell$. In our definition of congruences modulo~$\ell^n$
we are led by three requirements:
(1) If $n=1$, we want that $\alpha \equiv \beta \mod \ell$ if and only if the
reductions of $\alpha$ and $\beta$ are equal in~$\Fbar_\ell$.
(2) If $\alpha$ and $\beta$ are elements of some finite unramified extension $K/\QQ_\ell$,
then we want $\alpha \equiv \beta \mod \ell^n$ if and only of $\alpha - \beta \in (\pi_K^n)$.
(3) We want the definition to be independent of any choice of $K/\QQ_\ell$ containing
$\alpha$ and $\beta$.

We propose the following definition.

\begin{defi}\label{def:cong}
Let $n \in \NN$.
Let $\alpha,\beta \in \Zbar_\ell$. We say that $\alpha$ is congruent to $\beta$
modulo~$\ell^n$, for which we write $\alpha \equiv \beta \mod \ell^n$,
if and only if $v_\ell(\alpha - \beta) > n-1$.
\end{defi}

Note that this definition satisfies our three requirements.
Note also the trivial equivalence
\begin{equation}\label{eq:cong}
\alpha \equiv \beta \mod \ell^n \Leftrightarrow \lceil v_\ell(\beta - \alpha) \rceil \ge n.
\end{equation}
In the sequel of this article we will often speak of congruences
modulo~$\ell^n$ of (global) algebraic integers by fixing an
embedding $\Qbar \hookrightarrow \Qbar_\ell$. The same notation will
be used also in this situation without further comments.

\subsection{Interpretation in terms of ring extensions}

In this section we propose an interpretation of the above definition of congruences
modulo~$\ell^n$ in terms of ring extension of $\ZZ/\ell^n\ZZ$.
This interpretation gives us a much better algebraic handle for working with such
congruences because we will be able to use equality instead of congruence.
We were led to Definition~\ref{def:cong} by the following consideration:
Let $K/\QQ_\ell$ be a finite extension and $n \in \NN$.
What is the minimal~$m$ such that
the inclusion $\ZZ_\ell \hookrightarrow \cO_K$ induces an injection of
$\ZZ/\ell^n\ZZ$ into $\cO_K/(\pi_K^m)$?
In order to formulate the answer, we introduce a function.

\begin{defi}
Let $L/K/\QQ_\ell$ be finite field extensions and let $e_{L/K}$ denote the
ramification index of $L/K$. For $n \in \NN$, let $\gamma_{L/K} (n) = (n-1) e_{L/K} + 1$.
\end{defi}

This function satisfies the following simple properties:
\begin{enumerate}[(i)]
\itemsep=0cm plus 0pt minus 0pt
\parsep=0.1cm
\item\label{pi} For $n=1$, we have $\gamma_{L/K}(1) = 1$.
\item\label{pii} If $L/K$ is unramified, then $\gamma_{L/K}(n) = n$.
\item\label{piii} For extensions $M/L/K$, we have {\em multiplicativity}:
$\gamma_{M/K} (n) = \gamma_{M/L} (\gamma_{L/K}(n))$.
\item\label{piv} For extensions $L/K$, the integer $\gamma_{L/K}(n)$ is the minimal one
such that the embedding $\cO_K \hookrightarrow \cO_L$ induces an
injection $\cO_K / (\pi_K^n)\hookrightarrow \cO_L / (\pi_L^{\gamma_{L/K}(n)})$.
\item\label{pv} For $\alpha,\beta \in K/\QQ_\ell$ we have:
$$v_K(\alpha-\beta) \ge \gamma_{K/\QQ_\ell}(n) \Leftrightarrow v_\ell(\alpha-\beta) > n-1
\Leftrightarrow \alpha \equiv \beta \mod \ell^n.$$
\end{enumerate}

Note that (\ref{pi})--(\ref{piii}) precisely correspond to the requirements (1)--(3) from
Section~\ref{subsec:def}. By (\ref{piv}) we have produced {\em ring extensions}
$$  \ZZ/\ell^n \ZZ \hookrightarrow \cO_K / (\pi_K^{\gamma_{K/\QQ_\ell}(n)})
                    \hookrightarrow \cO_L / (\pi_L^{\gamma_{L/\QQ_\ell}(n)}).$$
Property~(\ref{pv}) immediately yields a reformulation of the congruence
of $\alpha$ and $\beta$ modulo~$\ell^n$
as an equality in the residue ring $\cO_K / (\pi_K^{\gamma_{K/\QQ_\ell}(n)})$.

In order to interpret congruences as equalities without always having to choose
some finite extension of~$\QQ_\ell$, we now make the following construction,
which for $n=1$ boils down to~$\Fbar_\ell$.
We define
$$\overline{\ZZ/\ell^n\ZZ} := \ilim{K} \cO_K/(\pi_K^{\gamma_{K/\QQ_\ell}(n)}),$$
where $K$ runs through all subextensions of $\Qbar_\ell$ of finite degree over~$\QQ_\ell$
and the inductive limit is taken with respect to the maps in~(\ref{piv}).
The natural projections
$\cO_K \twoheadrightarrow \cO_K/(\pi_K^{\gamma_{K/\QQ_\ell}(n)})$
give rise to a surjective ring homomorphism
$$ \pi_n: \Zbar_\ell \twoheadrightarrow \overline{\ZZ/\ell^n\ZZ}.$$
Now we can make another reformulation of our definition of congruences modulo~$\ell^n$:
Let $\alpha,\beta \in \Zbar_\ell$. Then we have
$$ \alpha \equiv \beta \mod \ell^n \Leftrightarrow \pi_n(\alpha) = \pi_n(\beta).$$
In the sequel, we will always choose the $\pi_n$ in a compatible way, i.e.\ if $m < n$ we want
$\pi_m$ to be the composition of $\pi_n$ with the natural map
$\overline{\ZZ/\ell^n\ZZ} \twoheadrightarrow \overline{\ZZ/\ell^m\ZZ}$.

\begin{rem}
We also point out a disadvantage of our choice of $\gamma_{K/\QQ_\ell}(n)$, namely that it is
not additive. This fact prevents us from defining a valuation on~$\Zbar_\ell$
by saying that the valuation of $a \in \Zbar$ is equal to the maximal
$n$ such that $\pi_n(a) = 0$. Defining $\gamma_{K/\QQ_\ell}(n)$ as
$n$ times the ramification index $e_{K/\QQ_\ell}$ would have avoided that problem.
But then $\gamma(1) = e_{K/\QQ_\ell} \neq 1$,
in general, which is not in accordance with the usual usage of modulo~$\ell$.
This other possibility can be understood as $\Zbar_\ell/\ell^n \Zbar_\ell$.
\end{rem}

\subsection{Computing congruences modulo $\ell^n$}

If one does not require one fixed embedding into the complex numbers, algebraic
integers are most easily represented by their minimal polynomials.
Thus, it is natural to study congruences between algebraic integers entirely
through their minimal polynomials. This is the point of view that we adapt and
it leads us to consider the following problem.

\begin{problem}\label{pproblem}
We fix, once and for all, for every~$n$ compatibly, ring homomorphisms
$\pi_n: \Zbar \hookrightarrow \Zbar_\ell \twoheadrightarrow \overline{\ZZ/\ell^n\ZZ}$.
Let $P,Q \in \ZZ[X]$ be two coprime monic polynomials and let $n \in \NN$.

How can we decide the validity of the following assertion?
\begin{quote}
``There exist $\alpha, \beta \in \Zbar$ such that
\begin{enumerate}[(i)]
\item $P(\alpha)=Q(\beta)=0$ and
\item $\pi_n(\alpha) = \pi_n(\beta)$ (i.e.\ $\alpha \equiv \beta \mod \ell^n$).''
\end{enumerate}
\end{quote}
\end{problem}

In this article, we will give two algorithms for treating this problem.
The first one arose from the idea that one could try to use greatest common divisors.
This notion seems to be the right one for $n=1$,
but it is not well behaved for $n > 1$ since the ring $\ZZ/\ell^n\ZZ[X]$
is not a principal ideal domain. However,
the algorithm for approximating greatest common divisors of two
polynomials over $\ZZ_\ell$ presented in Appendix~A of~\cite{FPR} led us to
consider the notion of {\em congruence number} or {\em reduced resultant}. It can
be used to give quite a fast algorithm, which, however, does not always give a complete
answer.

The second algorithm, which we call the {\em Newton polygon method}, always solves
Problem~\ref{pproblem} but tends to be slower (experimentally).
Its basic idea was suggested to us by Michael Stoll after a talk of the second author
and was immediately put into practice. However, since the first version of this article
had already been finished, the algorithm was not included in it, so that
it was again suggested to us by one of the referees.
In this section we will present both algorithms in detail.

It should be pointed out explicitly that Problem~\ref{pproblem} cannot be solved completely by
considering only the reductions of $P$ and $Q$ mod $\ell^n$ if $n>1$. This is a major
difference to the case $n=1$. The difference is due to the fact that in the problem
we want $\alpha$ and $\beta$ to be zeros of $P$ and~$Q$: if $\overline{\alpha}$ and
$\overline{\beta}$ are elements in $\overline{\ZZ/\ell^n\ZZ}$ such that inside that
ring $P(\overline{\alpha})=Q(\overline{\beta})=0$, then it is not clear if they are
reductions of zeros of $P$ and~$Q$.

\subsubsection*{Congruence number}

The congruence number of two integral polynomials provides an upper
bound for congruences in the sense of Problem~\ref{pproblem}. It is defined in
such a way that it can easily be calculated on a computer.

\begin{defi}
Let $R$ be any commutative ring. By $R[X]_{<n}$ we denote the $R$-module of
polynomials of degree less than~$n$.
Let $P,Q \in R[X]$ be two polynomials of degrees $m$ and~$n$, respectively.
The {\em Sylvester map} is the $R$-module homomorphism
$$ R[X]_{<n} \oplus R[X]_{<m} \to R[X]_{<(m+n)}, (r,s) \mapsto rP+sQ.$$
\end{defi}

If $R$ is a field, then the monic polynomial of smallest degree in the image
of the Sylvester map is the greatest common divisor of $P$ and~$Q$.
In particular, with $R$ a factorial integral domain and $P,Q$ primitive polynomials,
the Sylvester map is injective if and only if $P$ and $Q$ are coprime.
Consequently, if $P,Q \in \ZZ[X]$ are primitive coprime polynomials,
then any non-zero polynomial of smallest degree is a constant polynomial.

\begin{defi}
Let $P,Q \in \ZZ[X]$ be coprime polynomials.
We define the {\em congruence number}
$c(P,Q)$ of $P$ and $Q$ as the smallest positive integer $c$ such that
the constant polynomial~$c$ is in the image of the Sylvester map of $P$ and~$Q$.
\end{defi}

We remark that for monic coprime polynomials $P$ and $Q$
via polynomial division the principal ideal $(c(P,Q))$ can be seen to be
equal to the intersection of the ideal of constant integral polynomials with the
ideal in $\ZZ[X]$ generated by all polynomials
$rP+sQ$ when $r,s$ run through all of $\ZZ[X]$.
In \cite{Pohst} the congruence number is called the {\em reduced resultant}.
Note that in general the reduced resultant is a proper divisor of the resultant.
It makes sense to replace $\ZZ$ by $\ZZ_\ell$ everywhere and to define a congruence
number as a constant polynomial in the image of the Sylvester map having the lowest
$\ell$-adic valuation. Although this element is not unique, its valuation is.

The congruence number gives an upper bound for the $n$ in Problem~\ref{pproblem}:

\begin{prop}\label{ub}
Let $P,Q \in \ZZ[X]$ be coprime polynomials and let $\ell^n$ be the exact
power of~$\ell$ dividing $c(P,Q)$.
Then there are no $\alpha, \beta \in \Zbar$ such that
\begin{enumerate}[(i)]
\item $P(\alpha)=Q(\beta)=0$ and
\item $\pi_m(\alpha) = \pi_m(\beta)$ (i.e.\ $\alpha \equiv \beta \mod \ell^n$)
for any $m > n$.
\end{enumerate}
\end{prop}

\pf
By assumption there exist $r,s \in \ZZ[X]$ such that
$ c = c(P,Q) = rP + sQ$.
Let $\alpha, \beta \in \Zbar$ be zeros of~$P$ and $Q$, respectively, such that
$\pi_m(\alpha)=\pi_m(\beta)$. We obtain
$$\pi_m(c) = \pi_m\big(r(\alpha)P(\alpha) + s(\alpha)Q(\alpha)\big) \\
         = \pi_m\big(s(\alpha)\big) \pi_m\big(Q(\alpha)\big) \\
         = \pi_m\big(s(\beta)\big) \pi_m\big(Q(\beta)\big) = 0.$$
This means that $\ell^m$ divides~$c$, whence $m \le n$.
\qed

\subsubsection*{On the computation of the congruence number}

The idea for the computation of the congruence number is very simple:
we use basic linear algebra and the Sylvester matrix.
The point is that the Sylvester map is described by the standard Sylvester matrix~$S$
of $P$ and $Q$
(or rather its transpose if one works with column vectors) for the standard bases of
the polynomial rings.
We describe in words the straight forward algorithm
for computing the congruence number $c(P,Q)$
as well as for finding polynomials $r,s$ such that $c(P,Q) = rP+sQ$ with
$\deg(r) < \deg(Q)$ and $\deg(s) < \deg(P)$.
The algorithm consists of bringing $S$ into row echelon (or Hermite) form, i.e.\
one computes an invertible integral matrix~$B$ such that
$B S$ has no entries below the diagonal.
The congruence number $c(P,Q)$ is (the absolute value of) the bottom right entry of $BS$ and
the coefficients of $r$ and~$s$ are the entries in the bottom row of~$B$.
This algorithm works over the integers and over $\ell$-adic rings with a certain precision,
i.e.\ $\ZZ/\ell^n\ZZ$.

We note that by reducing $BS$ modulo~$\ell$, one can read off the greatest common
divisor of the reductions of $P$ and $Q$ modulo~$\ell$: its coefficients (up to normalization)
are the entries in the last non-zero row of the reduction of $BS$ modulo~$\ell$.
This has the following trivial, but noteworthy consequence.

\begin{cor}\label{modell}
Suppose that $P$ and $Q$ are primitive coprime polyomials in $\ZZ[X]$.
Then $P$ and $Q$ have a non-trivial common divisor modulo~$\ell$
if and only if the congruence number of $P$ and $Q$ is divisible by~$\ell$. \qed
\end{cor}

\subsubsection*{Applications of the congruence number}

We now examine when the congruence number is enough to solve Problem~\ref{pproblem}
for given $P,Q$ and for all~$n$.
In cases when it is not, we will give a lower bound for the maximum~$n$ for which
the assertions of the problem are satisfied.

We start with the observation that the congruence number suffices to solve
our problem for $n=1$.

\begin{prop}\label{neqone}
Let $n=1$. Assume that $P$ and $Q$ are coprime monic polynomials in $\ZZ[X]$.
The assertion in Problem~\ref{pproblem} is satisfied if and only if the congruence number
$c(P,Q)$ is divisible by~$\ell$.
\end{prop}

\pf
The calculations of the proof of Proposition~\ref{ub} show that if the assertion is
satisfied, then $\ell$ divides~$c(P,Q)$.
Conversely, if $\ell$ divides~$c(P,Q)$ then by Corollary~\ref{modell} the reductions
of $P$ and $Q$ have a non-trivial common divisor and thus a common zero in $\Fbar_\ell$.
All zeros in $\Fbar_\ell$ lift to zeros in $\Zbar_\ell$.
\qed
\medskip

We fix an embedding $\Qbar \hookrightarrow \Qbar_\ell$.
Our further treatment will be based on the following simple observation.
Let $M \subset \Qbar$ be any number field containing all the roots of the monic coprime
polynomials $P, Q \in \ZZ[X]$
and let $c = c(P,Q) = rP+sQ$ with $r,s \in \ZZ[X]$, $\deg(r) < \deg(Q)$, $\deg(s) < \deg(P)$ and
factor $Q(X) = \prod_i (X-\beta_i)$ in $\Zbar[X]$.
Then for $\alpha \in \Zbar$ such that $P(\alpha)=0$
we have
\begin{equation}\label{eqvM}
 v_M(c) = v_M\big(s(\alpha)\big) + \sum_i v_M(\alpha-\beta_i).
\end{equation}
Our aim now is to find a lower bound for the maximum of $v_M(\alpha-\beta_i)$
depending on $\pi_M(c)$.
For that we discuss the two summands in the equation separately.

We first treat $v_M\big(s(\alpha)\big)$.
By $\overline{F}$ we denote the reduction modulo~$\ell$ of an integral polynomial~$F$.

\begin{prop}\label{props}
Suppose that $\ell$ divides $c(P,Q)$.
\begin{enumerate}[(a)]
\item If $\overline{s}$ and $\overline{Q}$ are coprime, then $v_M\big(s(\alpha)\big) = 0$
for all $\alpha \in \Zbar$ with $\pi_1(Q(\alpha)) = 0$.

\item If one of $\overline{P}$ or $\overline{Q}$
does not have any multiple factors, then there is $\alpha \in \Zbar$ such that
$P(\alpha)= 0$, $\pi_1(Q(\alpha)) = 0$ and $v_M(s(\alpha)) = 0$, or
there is $\beta \in \Zbar$ such that
$Q(\beta)= 0$, $\pi_1(P(\beta)) = 0$ and $v_M(r(\beta)) = 0$.

\item If $\overline{P}$ is an irreducible polynomial in $\FF_\ell[X]$ and
$Q$ is irreducible in $\ZZ_\ell[X]$, then
$\overline{s}$ and $\overline{Q}$ are coprime and $v_M\big(s(\alpha)\big) = 0$
for all $\alpha \in \Zbar$ with $\pi_1(Q(\alpha)) = 0$.
\end{enumerate}
\end{prop}

\pf
(a) Since $\overline{s}$ and $\overline{Q}$ are coprime, the reduction of~$\alpha$
cannot be a root of both of them.

(b) We prove that there exists $y \in \Fbar_\ell$ which is
a common zero of $\overline{P}$ and $\overline{Q}$, but not
a common zero of $\overline{r}$ and $\overline{s}$ at the same time.
Assume the contrary, i.e.\ that $\overline{r}(y) = \overline{s}(y)=0$
for all $y \in \Fbar_\ell$ with $\overline{P}(y) = \overline{Q}(y)=0$.
Let $\overline{G} \in \FF_\ell[X]$ be the monic polynomial of smallest degree annihilating
all $y \in \Fbar_\ell$ with the property $\overline{P}(y) = \overline{Q}(y)=0$.
Then $\overline{G}$ divides $\overline{P}$, $\overline{Q}$ as well as by assumption
$\overline{r}$ and $\overline{s}$. Hence, we have
$$ 0 = \overline{r}\overline{P} + \overline{s}\overline{Q}
     = \overline{G}^2 \big( \overline{r_1}\overline{P_1} + \overline{s_1}\overline{Q_1}  \big)$$
with certain polynomials
$\overline{r_1}, \overline{P_1}, \overline{s_1}, \overline{Q_1} \in \FF_\ell[X]$.
We obtain the equation
\begin{equation}\label{eq:gcd}
 0 = \overline{r_1}\overline{P_1} + \overline{s_1}\overline{Q_1}
\end{equation}
and we also have $\deg(\overline{r_1}) < \deg(\overline{Q_1})$
and $\deg(\overline{s_1}) < \deg(\overline{P_1})$.
As either $\overline{P}$ or $\overline{Q}$ does not have any multiple factor,
it follows that $\overline{P_1}$ and $\overline{Q_1}$ are coprime.
This contradicts Equation~\ref{eq:gcd}.

Hence, we have $y \in \Fbar_\ell$ with $\overline{P}(y) = \overline{Q}(y)=0$
and $\overline{r}(y) \neq 0$ or $\overline{s}(y) \neq 0$.
If $\overline{r}(y) \neq 0$ then we lift $y$ to a zero $\beta$ of~$Q$.
In the other case we lift $y$ to a zero $\alpha$ of~$P$.

(c) The assumptions imply that $\overline{Q} = \overline{P}^a$ for some~$a$.
As the degree of $s$ is smaller than the degree of~$P$,
it follows that $\overline{s}$ and $\overline{P}$ are coprime.
Thus also, $\overline{s}$ and $\overline{Q}$ are coprime and we conclude by~(a).
\qed
\medskip

We now treat the term $\sum_i v_M(\alpha-\beta_i)$.

\begin{prop}\label{propg}
Suppose that $\ell$ divides $c(P,Q)$ and that $\alpha$ is a root of~$P$
which is congruent to some root of~$Q$ modulo~$\ell$ (which exists by Proposition~\ref{neqone}).
Assume without loss of generality that $\beta_1$ is a root of~$Q$ which is closest to~$\alpha$,
i.e.\ such that $v_M(\alpha-\beta_1) \ge v_M(\alpha-\beta_i)$ for all~$i$.

\begin{enumerate}[(a)]
\item Suppose that $\overline{Q}$ has no multiple factors (i.e.\ the discriminant of~$Q$
is not divisible by~$\ell$, or, equivalently, the congruence number of $Q$ and $Q'$
is not divisible by~$\ell$).

Then $\sum_i v_M(\alpha-\beta_i) = v_M(\alpha-\beta_1)$.

\item In general we have
$v_M(\alpha-\beta_1) \ge \lceil \frac{1}{\deg(Q)} \big( \sum_i v_M(\alpha-\beta_i)\big)  \rceil$.
\end{enumerate}
\end{prop}

\pf
(a) If $\overline{Q}$ does not have any multiple factors, then $v_M(\beta_1 - \beta_i) = 0$
for all $i \neq 1$. Consequently, $v_M(\alpha-\beta_i) = v_M(\alpha-\beta_1 + \beta_1-\beta_i)=0$
for $i \neq 1$.

(b) is trivial.
\qed
\medskip

We summarise of the preceding discussion in the following corollary,
solving Problem~\ref{pproblem} if $\overline{P}$ and $\overline{Q}$ do not have
any multiple factors, and giving a partial answer in the other cases.

\begin{cor}\label{cor-alg}
Let $P,Q$ be coprime monic polynomials in $\ZZ[X]$ (or $\ZZ_\ell[X]$) and
let $\ell^n$ be the highest power of~$\ell$ dividing the congruence number $c := c(P,Q)$ and let
$r,s \in \ZZ[X]$ (or $\ZZ_\ell[X]$) be polynomials such that
$c = rP+sQ$ with $\deg(r) < \deg(Q)$ and $\deg(s) < \deg(P)$.

\begin{enumerate}[(a)]
\item \label{part-exclude} If $n = 0$, then no root of $P$ is congruent modulo~$\ell$
to a root of~$Q$.
\item \label{part-one} If $n=1$, then there are $\alpha$, $\beta$ in $\Zbar$ (in $\Zbar_\ell$, respectively) with
$P(\alpha)=Q(\beta)=0$ such that they are congruent modulo~$\ell$,
and there are no $\alpha_1$, $\beta_1$ in $\Zbar$ (in $\Zbar_\ell$, respectively)
with $P(\alpha)=Q(\beta)=0$ such that they are congruent modulo~$\ell^2$.
\item \label{part-nice} Suppose now that $n \ge 1$ and that one of the
following properties holds:
\begin{enumerate}[(i)]
\item $\overline{P}$ does not have any multiple factors and
$\overline{Q}$ does not have any multiple factors (i.e.\ $\ell \nmid c(P,P')$ and
$\ell \nmid c(Q,Q')$).
\item $\overline{Q}$ does not have any multiple factors and $\overline{s}$
and $\overline{Q}$ are coprime.
\item $\overline{P}$ does not have any multiple factors and $\overline{r}$
and $\overline{P}$ are coprime.
\end{enumerate}
Then there are $\alpha$, $\beta$ in $\Zbar$ (in $\Zbar_\ell$, respectively) with
$P(\alpha)=Q(\beta)=0$ such that they are congruent modulo $\ell^n$ and
there are no $\alpha_1$, $\beta_1$ in $\Zbar$ (in $\Zbar_\ell$, respectively)
with $P(\alpha_1)=Q(\beta_1)=0$ such that they are congruent modulo~$\ell^{n+1}$.
\item\label{part-lb}
Suppose that $n \ge 1$.
\begin{enumerate}[(i)]
\item If $\overline{s}$ and $\overline{Q}$ are coprime,
let $m = \lceil \frac{n}{\deg(Q)} \rceil$.
\item If $\overline{r}$ and $\overline{P}$ are coprime,
let $m = \lceil \frac{n}{\deg(P)} \rceil$.
\item If (i) and (ii) do not hold, let $m=1$
\end{enumerate}
Then there are $\alpha$, $\beta$ in $\Zbar$ (in $\Zbar_\ell$, respectively) with
$P(\alpha)=Q(\beta)=0$ such that they are congruent modulo $\ell^m$ and
there are no $\alpha_1$, $\beta_1$ in $\Zbar$ (in $\Zbar_\ell$, respectively)
with $P(\alpha_1)=Q(\beta_1)=0$ such that they are congruent modulo~$\ell^{n+1}$.
\end{enumerate}
\end{cor}

\pf
In the proof we use the notation introduced above.
The upper bounds in (\ref{part-one})-(\ref{part-lb}) were proved in Proposition~\ref{ub}.

(a) follows from Proposition~\ref{neqone}.

(b) The existence of a congruence follows from Corollary~\ref{modell}.

(c) In case (i), by Proposition~\ref{props}~(b)
we can choose $\alpha, \beta \in \Zbar$ congruent modulo~$\ell$ with $P(\alpha)=0$
and $\beta \in \Zbar$ with $Q(\beta)=0$ such that $v_M(s(\alpha)) =0$ or
$v_M(r(\beta))=0$. Without loss of generality
(after possibly exchanging the roles of $(P,r)$ and $(Q,s)$)
we may assume the former case.
In case (ii), by Proposition~\ref{props}~(a) any $\alpha \in \Zbar$
with $P(\alpha)=0$ and $\pi_1(Q(\alpha))=0$ will satisfy $v_m(s(\alpha))=0$.
In both cases, from Proposition~\ref{propg}
and Equation~\ref{eqvM} we obtain the equality
$$v_M(c)=v_M(\ell^n) = v_M(\alpha - \beta_1),$$
where $\beta_1$ comes from Proposition~\ref{propg}.
This gives the desired result.
Case (iii) is just case (ii) with the roles of $(P,r)$ and $(Q,s)$ interchanged.

(d) also follows from Propositions \ref{props} and~\ref{propg} and Equation~\ref{eqvM}.
More precisely, in case (i) we have the inequality
$$v_M(\alpha-\beta_1) \ge \lceil \frac{v_M(c)}{\deg(Q)} \rceil
= \lceil \frac{e n}{\deg(Q)} \rceil
\ge \big(\lceil \frac{n}{\deg(Q)} \rceil -1  \big) e + 1
= \gamma_{M/\QQ_\ell} (\lceil \frac{n}{\deg(Q)} \rceil),$$
where $e$ is the ramification index of $M/\QQ_\ell$. Hence,
$\pi_m(\alpha - \beta_1) = 0$ with $m = \lceil \frac{n}{\deg(Q)} \rceil$.
Case (ii) is case (i) with the roles of $(P,r)$ and $(Q,s)$ interchanged.
\qed

\begin{rem}
It is straightforward to turn Corollary~\ref{cor-alg} into an algorithm.
Say, $P,Q \in \ZZ[X]$ are coprime monic polynomials.
First we compute the congruence numbers $c(P,P')$ and $c(Q,Q')$. If any of these is zero,
then we factor $P$ (respectively, $Q$) in $\ZZ[X]$ into irreducible polynomials
$P = \prod_i P_i$ (respectively, $Q = \prod_j Q_j$). We then treat any pair $(P_i,Q_j)$
separately and return the maximum upper and the maximum lower bound for congruences of zeros.
For simplicity of notation, we now call the pair $(P,Q)$.

Now we compute the congruence numbers $c = c(P,Q)$ and $c_P =
c(P,P')$ as well as $c_Q = c(Q,Q')$, all of which are non-zero by
assumption. Along the way we also compute polynomials $r,s \in
\ZZ[X]$ such that $c = rP+sQ$ and $\deg(r) < \deg(Q)$ and $\deg(s) < \deg(P)$.
For each prime power $\ell^n$ (with $n \ge 1$) exactly dividing~$c$ we do the following.
If $\ell$ does not divide $c_P c_Q$, then we are in case~(\ref{part-nice})(i) and we
know that there are $\alpha, \beta \in \Zbar$ such that
$P(\alpha)=0=Q(\beta)$ and $\pi_n(\alpha) = \pi_n(\beta)$.
This is best possible and we have obtained a complete answer to Problem~\ref{pproblem}.
If $\ell$ is coprime to $c_P$ or $c_Q$, we check whether we are
in case (\ref{part-nice})(ii) or (\ref{part-nice})(iii). Then we also
obtain equality of the upper and lower bound and thus
a complete answer to Problem~\ref{pproblem}.
If we are in neither of these cases, then we use the much weaker lower bounds
of part~(\ref{part-lb}). In order to get a best possible result in this case, too,
one can make use of the Newton polygon method to be described next.
\end{rem}

\subsubsection*{Newton polygon method}

We now present the second algorithm for treating Problem~\ref{pproblem}.
The basic idea of this algorithm was suggested to us by Michael Stoll.
Let still $P,Q \in \ZZ[X]$ be coprime monic polynomials.
Consider factorisations in $\Zbar[X]$:
$$P(X) = \prod_{i=1}^u(X-\alpha_i) \textnormal{ and } Q(X) = \prod_{j=1}^v(X-\beta_j).$$
Now take $Q(X+Y) = \prod_{j=1}^v(X-(\beta_j-Y))$, considered as a polynomial
in $X$ with coefficients in~$\ZZ[Y]$ and let $F(Y)$ be the resultant of $P(X)$ and $Q(X+Y)$
with respect to the variable~$X$. By well known properties of the resultant one has
$$ F(Y) = \pm \prod_{i=1}^u\prod_{j=1}^v (Y-(\beta_j-\alpha_i)).$$
Hence, the roots of $F(Y)$ are precisely the differences of the roots of $P$ and~$Q$.
Thus, the slopes of the Newton Polygon of $F(Y) \in \ZZ_\ell[Y]$
are the $v_\ell(\beta_j-\alpha_i)$.
We obtain the following result, solving Problem~\ref{pproblem}.

\begin{prop}
Let $P,Q \in \ZZ[X]$ be coprime monic polynomials and set
$n := \lceil s \rceil$, where $s$ is the biggest slope of the Newton polygon of
the polynomial $F \in \ZZ_\ell[Y]$ defined above.

Then there are $\alpha, \beta \in \Zbar$ such that
\begin{enumerate}[(i)]
\item $P(\alpha)=Q(\beta)=0$ and
\item $\pi_n(\alpha) = \pi_n(\beta)$ (i.e.\ $\alpha \equiv \beta \mod \ell^n$).
\end{enumerate}
Moreover, $n$ is the biggest integer satisfying this property.
\end{prop}

\pf
Let $\alpha, \beta \in \Zbar$ with $P(\alpha)=Q(\beta)=0$ such that the slope
of $\beta - \alpha$ is equal to~$s$, i.e.\ $v_\ell(\beta-\alpha)=s$
(subject to the fixed embedding $\Qbar \hookrightarrow \Qbar_\ell$).
The proposition is an immediate consequence of Definition~\ref{def:cong}
and Equation~\ref{eq:cong}.
\qed

\section{Modular forms and Galois representations modulo $\ell^n$}\label{sec-modular}

In this section, we apply the methods from Section~\ref{sec-cong} to the study
of congruences of modular forms and modular Galois representations modulo~$\ell^n$.

As in Section~\ref{sec-cong}, we keep ring homomorphisms
$\pi_n: \Zbar \hookrightarrow \Zbar_\ell \twoheadrightarrow \overline{(\ZZ/\ell^n\ZZ)}$,
compatibly for~$n$, fixed.
In this section, we restrict to $\Gamma_0(N)$ for simplicity. Everything can be generalised
without any problems to $\Gamma_1(N)$ with the obvious modifications. Moreover,
also for the simplicity of the exposition all our modular forms are cusp forms.

\subsection{Modular forms modulo $\ell^n$}

For studying the notion of congruences modulo~$\ell^n$ of modular forms it is useful
to introduce the terminology of modular forms over~$\ZZ/\ell^n \ZZ$ or, in abuse of language,
modular forms modulo~$\ell^n$.
In contrast to the case $n=1$, one must be aware that lifting of modular forms over $\ZZ/\ell^n\ZZ$ to characteristic zero is not automatic. This will be reflected in our notions.
We let $S_k(\Gamma_0(N))$ denote the $\CC$-vector space of holomorphic cuspidal
modular forms of weight~$k$ and level~$N$.

\begin{defi}\label{defi-eigenform}
Let $\TT := \TT_k(\Gamma_0(N))$
 be the $\ZZ$-subalgebra of $\End_\CC(S_k(\Gamma_0(N)))$ generated by all
the Hecke operators $T_n$, $n \in \NN$.
\begin{enumerate}[(i)]
\item A {\em modular form of weight~$k$ and level~$N$ over $\ZZ/\ell^n\ZZ$
 (or modulo $\ell^n$)} is a $\ZZ$-module homomorphism
$ f: \TT \to \overline{(\ZZ/\ell^n\ZZ)}$.
\item A modular form $f$ over $\ZZ/\ell^n\ZZ$ is a {\em weak Hecke eigenform} if $f$
is a ring homomorphism.
\item \label{part:strong} A weak Hecke eigenform~$f$ over $\ZZ/\ell^n\ZZ$
is a {\em strong Hecke eigenform} if $f$ factors into ring homomorphisms
$\TT \to \Zbar_\ell \xrightarrow{\pi_n} \overline{(\ZZ/\ell^n\ZZ)}$.
\item Any normalised holomorphic Hecke eigenform $f = q + \sum_{m \ge 2} a_m(f) q^m$
(with $q = e^{2\pi i z}$ and $a_m \in \Zbar$)
gives rise to a strong Hecke eigenform over $\ZZ/\ell^n\ZZ$ via
$\TT \xrightarrow{T_m \mapsto a_m} \Zbar \xrightarrow{\pi_n} \overline{(\ZZ/\ell^n\ZZ)}$.
This modular form will be referred to as the {\em reduction of $f$ modulo~$\ell^n$}.
\item If the reductions modulo~$\ell^n$ of two normalised holomorphic eigenforms $f$ and $g$ agree,
then we say that $f$ and $g$ are {\em congruent modulo~$\ell^n$}.
This is the same as the congruence $a_m(f) \equiv a_m(g) \mod~\ell^n$ for all~$m\in \NN$
with the notion of congruence from Section~\ref{sec-cong}.
If the congruence $a_p(f) \equiv a_p(g) \mod~\ell^n$ holds for all primes~$p$ but possibly
finitely many, we say that $f$ and $g$ are {\em congruent modulo~$\ell^n$ at almost all primes}.

\end{enumerate}
\end{defi}

\begin{rem}\label{rem-eigenform}
\begin{enumerate}[(a)]
\item It is often useful to think of a modular form~$f$ over $\ZZ/\ell^n\ZZ$
as the {\em $q$-expansion} $\sum_{n=1}^\infty f(T_n) q^n \in \overline{\ZZ/\ell^n\ZZ}[[q]]$.

\item As $\TT$ is a finitely generated (and free) $\ZZ$-module, every weak eigenform~$f$
can be factored as
$\TT \to \cO_K/(\pi_K^{\gamma_{K/\QQ_\ell}(n)}) \to \overline{\ZZ/\ell^n\ZZ}$
for a suitable $\ell$-adic field~$K$.

\item Let $f: \TT\xrightarrow{\phi} \Zbar_\ell \xrightarrow{\pi_n} \overline{\ZZ/\ell^n\ZZ}$
be a strong Hecke eigenform modulo~$\ell^n$.
The kernel of~$\phi$ is a minimal prime ideal $\fp$ of~$\TT$.
As such, it corresponds to a $\Gal(\Qbar/\QQ)$-conjugacy class of holomorphic
Hecke eigenforms, since
$L := \Frac(\TT/\fp) \subseteq \Qbar$ is a number field (recall that
$\TT$ is a finitely generated free $\ZZ$-module)
and $\fp$ is the kernel of the ring homomorphism
$$\TT \twoheadrightarrow \TT/\fp \subset L \hookrightarrow \Qbar
\subset \CC,\;\; T_m \mapsto a_m,$$
which corresponds to the
normalised holomorphic eigenform $\sum_{m\ge 1} a_m e^{2\pi i m z}$ and depends on
the choice of the embedding $L \hookrightarrow \Qbar$.
Hence, the notion of strong Hecke eigenform modulo~$\ell^n$ implies that the form~$f$
is the reduction of a holomorphic Hecke eigenform modulo~$\ell^n$.

\item For $n=1$, the notion of weak and strong Hecke eigenform agree.
The reason is that the kernel of $f: \TT \to \Fbar_\ell$
is a maximal ideal, since the image of $f$ is a (finite) field. Every maximal ideal of~$\TT$
contains a minimal prime ideal~$\fp$ and, hence,
$f$ factors as
$\TT \to \TT/\fp \hookrightarrow \Zbar \hookrightarrow \Zbar_\ell \twoheadrightarrow \Fbar_\ell$.

\item Weak Hecke eigenforms need not be strong Hecke eigenforms in general.
See, for instance, Section~\ref{sec-weakstrong}.

\item Let $R$ be any ring. Since $\Hom_\ZZ(\TT,\ZZ) \otimes_\ZZ R \cong
\Hom_\ZZ(\TT,R)$ due to the freeness of $\TT$ as a finitely generated
$\ZZ$-module and since $\Hom_\ZZ(\TT,\ZZ)$ can be identified with the holomorphic
modular forms having integral Fourier expansions, any homomorphism
$f: \TT \to R$ (e.g.\ weak/strong eigenform) can be seen as an $R$-linear
combination of holomorphic modular forms (which are not necessarily eigenforms).

\item \label{rem-bad-galois} Another issue concerns
the absence of a good Galois theory for the extensions of
$\ZZ/\ell^n\ZZ$ discussed in Section~\ref{sec-cong}: Let $K$ be an $\ell$-adic field.
Not every ring homomorphism
$\cO_K \to \cO_K/(\pi_K^m)$ comes from a field homomorphism $K \to K$.
Suppose, for example, that $\cO_K = \ZZ_\ell[X]/(P(X))$ is the ring of integers
of a ramified extension of~$\QQ_\ell$.
If $\alpha$ is a root of~$P$ and if $m$ is big enough, then $\alpha+\pi^{m-1}$ is not
a root of~$P$, but nevertheless $P(\alpha+\pi^{m-1}) \in (\pi_K^m)$, whence sending
$\alpha$ to $\alpha+\pi^{m-1}$ uniquely defines a ring homomorphism $\cO_K \to \cO_K/(\pi_K^m)$,
which does not lift to a field automorphism $K\to K$.
Hence, a strong Hecke eigenform modulo~$\ell^n$ can give rise to many weak Hecke
eigenforms modulo~$\ell^n$.

\item Finally, we would like to point out a connection, as suggested by one of the
referees, between the congruence number and the congruence exponent of
modular abelian varieties defined in the paper~\cite{ARS} by Agashe, Ribet
and Stein and our notions.

Let $J$ be the Jacobian (over~$\QQ$) of some modular curve (say, $X_0(N)$)
and $A$, $B$ abelian subvarieties of~$J$ such that $J=A+B$ and $A \cap B$ is finite.
For the moment, let $\TT$ be the Hecke algebra of~$J$, i.e.\ the subring
of the endomorphism ring of~$J$ generated by all Hecke operators.
Denote by $\TT_A$ and $\TT_B$ the Hecke algebras of $A$ and~$B$, respectively.
The natural map $\phi: \TT \to \TT_A \oplus \TT_B$ given by sending an operator~$T$
to its restrictions to~$A$ and~$B$ is injective due to the
condition $J=A+B$. Thus, we can view~$\TT$
as an abelian subgroup of~$\TT_A \oplus \TT_B$, which has finite index, since $A \cap B$
is finite.
Agashe, Ribet and Stein define the {\em congruence exponent} (and the
{\em congruence number}) of~$A$
as the exponent (the number of elements) of the abelian group $(\TT_A \oplus \TT_B) /\TT$.
Note that the definition also depends on~$B$.

Now we establish the connection to our set-up. The Hecke algebra~$\TT$ is known to
be isomorphic to the Hecke algebra $\TT_2(\Gamma_0(N))$. Applying the functor
$\Hom_\ZZ(\cdot, \overline{\ZZ/\ell^n\ZZ})$, we obtain the exact sequence
\begin{multline*}
0 \to \Hom_\ZZ((\TT_A \oplus \TT_B)/\TT, \overline{\ZZ/\ell^n\ZZ}) \xrightarrow{\alpha}
\Hom_\ZZ(\TT_A, \overline{\ZZ/\ell^n\ZZ}) \oplus \Hom_\ZZ(\TT_B, \overline{\ZZ/\ell^n\ZZ})\\
\xrightarrow{\beta} \Hom_\ZZ(\TT,\overline{\ZZ/\ell^n\ZZ}).
\end{multline*}
Note that the term on the right is precisely the group of weight~$2$ modular forms
modulo~$\ell^n$ on $\Gamma_0(N)$ in our definition.
Let us now take two normalised
newforms $f$ and~$g$ in $S_2(\Gamma_0(N))$ in distinct Galois conjugacy classes
such that $f$ corresponds to a ring homomorphism $f:\TT_A \to \CC$ and
$g$ to $g:\TT_B \to \CC$. This is the case, for instance, if $A = (J/I_f J)^\vee$ and
$B = I_f J$, where $I_f$ is the kernel of the ring homomorphism $\TT \to \CC$
belonging to~$f$.
Assume that $f$ and~$g$ are congruent modulo~$\ell^n$. This means by definition that
$(f,-g)$ is in the kernel of~$\beta$. We analyse the element
$\psi \in \Hom_\ZZ((\TT_A \oplus \TT_B)/\TT, \overline{\ZZ/\ell^n\ZZ})$ such
that $\alpha(\psi) = (f,-g)$. It satisfies $\psi( (T_1,0) + \TT ) = f(T_1) - g(0) = 1$,
since $f$ is normalised. Consequently, $\ZZ/\ell^n\ZZ$ is in the image of~$\psi$.
Hence, $(\TT_A \oplus \TT_B)/\TT$ contains an element of order~$\ell^n$.
We conclude that $\ell^n$ divides the congruence exponent of~$A$ (and, of course,
also the congruence number).

\end{enumerate}
\end{rem}

\subsection{Galois Representations modulo $\ell^n$}\label{sec-rep}

We are interested in congruences modulo $\ell^n$ (in the sense of Section~\ref{sec-cong})
between $2$-dimensional $\ell$-adic Galois representations ($i=1,2$)
$$ \rho_i: \Gal(\Qbar/\QQ) \to \GL_2(\cO_{K_i}),$$
i.e.\ $\cO_{K_i}$ is the ring of integers of an $\ell$-adic field.
For that let $K$ be an $\ell$-adic field containing $K_1$ and~$K_2$.
We study the reductions of the representations modulo~$\ell^n$:
$$ \rhobar_i^{(n)}: \Gal(\Qbar/\QQ) \to \GL_2(\cO_K)\xrightarrow{\textnormal{nat.\ proj.}} \GL_2(\cO_K/(\pi_K^{\gamma_{K/\QQ_\ell}(n)})).$$

\begin{defi}
The representations $\rho_1$ and $\rho_2$ are called {\em congruent
modulo~$\ell^n$} if $\rhobar_1^{(n)}$ and $\rhobar_2^{(n)}$ are
isomorphic as
$(\cO_K/(\pi_K^{\gamma_{K/\QQ_\ell}(n)}))[\Gal(\Qbar/\QQ)]$-modules.
\end{defi}

\begin{rem}
The insistence on taking the natural projection is owed to the fact that there
may be `too many' maps from $\cO_K \to \cO_K/(\pi_K^{\gamma_{K/\QQ_\ell}(n)})$, as mentioned
in Remark~\ref{rem-eigenform}~(\ref{rem-bad-galois}).
\end{rem}

\begin{thm} \label{thm-galrep}
If the $\rho_i$ are residually absolutely irreducible, then they are congruent modulo~$\ell^n$
if and only if the traces of Frobenius elements agree, i.e.\
$\Tr(\rhobar_1^{(n)}(\Frob_p)) = \Tr(\rhobar_2^{(n)}(\Frob_p))$,
at a dense set of primes~$p$.
\end{thm}

\pf
Chebotarev's Theorem applied to the Proposition in \cite{M}, p.~253.
\qed
\medskip

Subject to a fixed choice $\Qbar \hookrightarrow \Qbar_\ell$, to a normalised
holomorphic eigenform~$f = \sum a_m q^m \in S_k(\Gamma_0(N))$
one can attach an $\ell$-adic Galois
representation $\rho_{f,\ell}: \Gal(\Qbar/\QQ) \to \GL_2(K)$
with some (suitably large) $\ell$-adic field~$K$. This Galois representation has
the properties that it is unramified outside $\ell$ and the level of~$f$ and the
trace of $\Frob_p$ is equal to~$a_p$ at all unramified primes~$p$.

\begin{prop}
Any weak or strong Hecke eigenform~$f: \TT \to \cO_K/(\pi_K^{\gamma_{K/\QQ_\ell}(n)})$
of level $N$ and weight~$k$ has an attached residual Galois representation $\rhobar_{f,\ell}$.
If $\rhobar_{f,\ell}$ is absolutely irreducible, $f$ gives
rise to a Galois representation
$$\rhobar_{f,\ell}^{(n)}: \Gal(\Qbar/\QQ) \to \GL_2(\cO_K/(\pi_K^{\gamma_{K/\QQ_\ell}(n)}))$$
which is unramified outside $\ell N$ and satisfies for every $p\nmid\ell N$
$$  \Tr(\rhobar_{f,\ell}^{(n)}(\Frob_p))=a_p,\text{ and }
   \Det(\rhobar_{f,\ell}^{(n)}(\Frob_p))=p^{k-1}, $$
where we write $a_p$ for the $p$-th coefficient of~$f$, i.e.\ $a_p=f(T_p)$.
\end{prop}

\pf Any weak modular form modulo~$\ell^n$ gives rise to a strong modular form
modulo~$\ell$ by reduction, and hence we dispose of $\rhobar_{f,\ell}$. If the
residual representation is absolutely irreducible,
Theorem~3 (p.~225) from~\cite{C} implies the existence of a
Galois representation
$$ \rho: \Gal(\Qbar/\QQ) \to \GL_2(\TT \otimes_\ZZ \ZZ_\ell)$$
with the desired properties. Note that $f$ factors as
$\TT \to \TT \otimes_\ZZ \ZZ_\ell \xrightarrow{f_1} \cO_K/(\pi_K^{\gamma_{K/\QQ_\ell}(n)})$.
It hence suffices to compose $\rho$ with the natural map coming from~$f_1$.
\qed

\subsection{Sturm bound modulo $\ell^n$}

If two Galois representations $\rhobar_i^{(n)}$ ($i=1,2$) as in the previous subsection come
from weak or strong modular forms modulo~$\ell^n$,
then one can decide whether they are equivalent
by comparing only finitely many coefficients, since one disposes of an effective
bound for the two modular forms modulo~$\ell^n$ to be equal. Such a bound is
given by the Sturm bound (\cite{Sturm}).

\begin{thm}\label{thm-fin-gen}
Let $\Gamma$ be a congruence group containing $\Gamma_1(N)$, let $k\ge 1$
and let $B$ be the \emph{Sturm bound} defined by
\[
B:=\frac{kb}{12}-\frac{b-1}{N},
\]
where $b=[\SL_2(\ZZ):\Gamma]$. The Hecke algebra $\TT$ acting on the
space $S_k(\Gamma)$ is generated as a $\ZZ$-module
by the Hecke operators $T_n$ for $1\le n\le B$.
Moreover, for $\Gamma = \Gamma_0(N)$ the algebra $\TT$ is generated as a $\ZZ$-algebra
by the~$T_p$ for the primes $p\le B$.
\end{thm}
\pf Theorem~$9.23$ and Remark $9.24$ from \cite{S}. \qed

\begin{thm}\label{thm-congr}
Let $f,g : \TT \to \cO_K/(\pi_K^{\gamma_{K/\QQ_\ell}(n)})$ be two weak or strong
Hecke eigenforms modulo~$\ell^n$ on $\Gamma_0(N)$ for some weight~$k$.
Let $b=[SL_2(\ZZ):\Gamma_0(N)]$.
If for all primes
\[
p\le\frac{k b}{12}-\frac{b-1}{N}
\]
we have
\[
f(T_p) = g(T_p) \;\;\;\;\;\;\;\;\;\text{(i.e.\ ``$a_p(f) \equiv a_p(g) \mod \ell^n$'')},
\]
then $f$ is equal to $g$ as a Hecke eigenform modulo~$\ell^n$.
\end{thm}

\pf 
As for $\Gamma = \Gamma_0(N)$ we have that $\TT$ is generated as a $\ZZ$-algebra by the Hecke
operators $T_p$ for the primes $p\le B$
(Theorem~\ref{thm-fin-gen}), it follows that $f$ and $g$ are uniquely determined
by their values at $T_p$ for primes $p \le B$.
\qed

\begin{rem}
The Sturm bound can easily be extended to modular forms with nebentype, see e.g.\
\cite{S}, Corollary~9.20.
\end{rem}

We mention that in \cite{CKR}, the Sturm bound is proved by other
means and is also extended to the situation when the two modular
forms have different weights. It is also useful to remark that the
Sturm bound for modular forms modulo~$\ell^n$ is also a direct
consequence of the Sturm bound for modular forms over~$\FF_\ell$ and
Nakayama's Lemma: If $\TT \otimes_\ZZ \FF_\ell$ is generated as
$\FF_\ell$-vector space by the Hecke operators $T_1,\dots,T_B$, then
$\TT \otimes_\ZZ \ZZ/\ell^n\ZZ$ is generated as a
$\ZZ/\ell^n\ZZ$-modulo by $T_1,\dots,T_B$, too.

\subsection{Application of degeneracy maps}

Theorem~\ref{thm-congr} gives a criterium for the Galois
representations attached to two Hecke eigenforms $f \in S_k(\Gamma_0(N))$ and
$g \in S_k(\Gamma_0(Nm))$ to be congruent modulo $\ell^n$
(under the assumption that the representations are residually irreducible).
However, most of the time when studying congruences of Galois representations
attached to modular forms $f$ and $g$, the assumptions of
Theorem~\ref{thm-congr} will not be fulfilled, as $f$ and $g$ will
typically differ at some prime dividing one of the levels. Hence, we
now propose a stronger criterion. In order to formulate it, we
introduce some straightforward notation.

\begin{defi}
Let $R$ be a commutative ring (in the sequel, either $R=\CC$, $R=\ZZ$
or $R$ is an extension of $\ZZ/\ell^n\ZZ$ as in Section~\ref{sec-cong})
and $d \in \NN$.
Let $N,m,n \in \NN$.
The {\em degeneracy map} for a positive divisor $d$ of~$m$ is defined to be the map
$$\phi_d: \Hom_\ZZ(\TT_k(\Gamma_0(N)),R) \to \Hom_\ZZ(\TT_k(\Gamma_0(Nm)),R)$$
which sends $f \in \Hom_\ZZ(\TT_k(\Gamma_0(N)),R)$ to the
homomorphism in $\Hom_\ZZ(\TT_k(\Gamma_0(Nm)),R)$ that maps $T_n$
to $\phi(T_{n/d})$, if $d$ divides~$n$, and to~$0$ otherwise.

Let $f:\TT_k(\Gamma_0(N)) \to R$ be a modular form over~$R$.
The {\em old space of $f$ over~$R$ in level~$Nm$} is defined as the
$R$-span of the image of $f$ under the degeneracy maps for each positive~$d\mid m$
inside $\Hom_\ZZ(\TT_k(\Gamma_0(Nm)),R)$.
\end{defi}

On $q$-expansions, the degeneracy map for~$d$ corresponds to the $R$-module endomorphism
of $R [[q]]$ given by $q \mapsto q^d$.
The degeneracy map~$\phi_d$ is well defined with $R=\ZZ$ by the classical theory of
modular forms (via the identification of $\Hom_\ZZ(\TT_k(\Gamma_0(N)),\ZZ)$
with those holomorphic cusp forms in $S_k(\Gamma_0(N))$ having integral
Fourier expansions) and due to the isomorphism
$\Hom_\ZZ(\TT_k(\Gamma_0(N)),\ZZ) \otimes_\ZZ R \cong \Hom_\ZZ(\TT_k(\Gamma_0(N)),R)$
it is well defined for all rings~$R$.

\begin{prop}\label{prop-use-sturm}
Let $f$ and $g$ be weak Hecke eigenforms modulo~$\ell^n$ of weight~$k$ for $\Gamma_0(N)$ and
$\Gamma_0(Nm)$, respectively, and assume that their residual Galois representations
are absolutely irreducible.

Then the Galois representations modulo~$\ell^n$ attached to $f$ and~$g$
are isomorphic if there is a weak Hecke eigenform $\tilde{f}$ modulo~$\ell^n$
in the oldspace of~$f$ modulo~$\ell^n$ in level~$Nm$ such that
$g(T_p) = \tilde{f}(T_p)$ (i.e.\ ``$a_p(g) \equiv a_p(\tilde{g}) \mod \ell^n$'')
for the primes~$p$ up to the Sturm bound for weight~$k$ and $\Gamma_0(Nm)$.
\end{prop}

\pf
The assumptions imply that the equality
$g(T_p) = f(T_p)$ holds for all primes~$p$ except possibly those
with $p$ dividing~$m$. Hence, we can conclude by
Theorem~\ref{thm-galrep}.
\qed
\medskip

Proposition~\ref{prop-use-sturm} gives rise to a
straightforward algorithm (see Section~\ref{sec-alg}), since the characteristic polynomials
of the Hecke operators at~$p \mid m$ on the oldspace of~$f$ can be described explicitly
as follows.
Let $f\in S_k(\Gamma_0(N))$ and $g \in S_k(\Gamma_0(Nm))$
be Hecke eigenforms. Suppose that $r$ is the maximum exponent
such that $p^r\mid m$.
Then $T_p$ acts on the old space of~$f$ in level $p^r N$ as the $(r+1)\times (r+1)$ matrix
\begin{equation}\label{bigmatrix}
\tilde{T}_p = \left(
\begin{array}{cccccc}
   a_p(f)    & 1      & 0 & 0 & \ldots & 0      \\
   -\delta p^{k-1} & 0      & 1 & 0 & \ldots & 0      \\
   0         & 0      & 0 & 1 & \ldots & 0      \\
   \vdots    &        &   &   &        & \vdots \\
   0         & \ldots & 0 & 0 & 0      & 1      \\
   0         & \ldots & 0 & 0 & 0      & 0      \\
\end{array}
\right)
\end{equation}
where $\delta=0$ if $p\mid N$ and $\delta=1$ otherwise (see \cite{W}).

Let $[f]$ be the $\ZZ$-span of the $\Gal(\Qbar/\QQ)$-conjugacy class of~$f$;
say that its rank is~$d$.
The operator $T_p$ acts on the image of~$[f]$ in level $mN$
as the $d\cdot(r+1)\times d\cdot(r+1)$ matrix resulting
from (\ref{bigmatrix}), in which we substitute every $0$ by the
$d\times d$ dimensional $0_d$ matrix, $1$ becomes the $d$-identity
$1_d$, the entry $a_p(f)$ is replaced by the $d\times d$ matrix of the Hecke operator
$T_p$ on~$[f]$, and $\delta$ is either $0_d$ or $1_d$.
Since all the elements below the diagonal are~$0$ for all the blocks
under the second line of blocks, we know that the
characteristic polynomial of this big matrix will be the product of
$X^{d(r-1)}$ and the characteristic polynomial of the block matrix
\begin{equation}\label{mittlematrix}
\left(
\begin{array}{c|c}
   T_p    & 1_d    \\
\hline
   -\delta p^{k-1}\cdot1_d & 0_d
\end{array}
\right).
\end{equation}
We now compute the characteristic polynomial of (\ref{mittlematrix}).
Let $P_{f,p}=\sum_{i=0}^dc_iX^i = \prod_{j=1}^d (X-a_j)$ be the characteristic polynomial
of the upper left block, where the $a_j$ lie in some algebraic closure.
With two polynomial variables~$\tilde{X},\tilde{Y}$ we hence have
$\prod_j (\tilde{X}-a_j\tilde{Y}) = \sum_i c_i \tilde{X}^i \tilde{Y}^{d-i}$.
We now plug in $\tilde{X} = X^2+\delta p^{k-1}$ and $\tilde{Y}=X$ and obtain
$$ \prod_{j=1}^d(X^2 - a_jX +\delta p^{k-1}) = 
\sum_{i=0}^d\bigg{(}c_iX^{d-i}(X^2+\delta p^{k-1})^i\bigg{)}. $$
By taking the Jordan normal form (over an algebraic closure) and rearranging the
matrix, we see that this is the characteristic polynomial of (\ref{mittlematrix}).
Hence, the characteristic polynomial $\tilde{P}_{f,p}$ of
\ref{bigmatrix} is
\begin{equation}\label{charpoly}
\tilde{P}_{f,p}=\sum_{i=0}^d\bigg{(}c_iX^{dr-i}(X^2+\delta p^{k-1})^i\bigg{)},
\end{equation}
which can be computed very quickly from $P_{f,p}$. Let us remark
that, if $p\mid N$, this polynomial is simply $X^{dr}\cdot P_{f,p}$
and, if $p\nmid N$ and $d=1$, then $\tilde{P}_{f,p}$ is $X^{r-1}$ times the
characteristic polynomial of the $p$-Frobenius element.

\begin{rem}
\begin{enumerate}[(a)]
\item It appears worthwhile to investigate the existence of
a partial converse to Proposition~\ref{prop-use-sturm}. A true
converse cannot hold if $f$ is in the lowest possible level,
since it is easy to construct a counter example if $n=1$, $k=2$ and
$\ell=2$ and there is a weight-$1$ form embedded into weight~$2$.
Under certain conditions (e.g.\ $k< \ell$ and $\ell \nmid Nm$) a
converse could conceivably exist.

To illustrate the problem with a particular example, let us consider the unique
Hecke eigenform~$f$ modulo~$2$ in level $\Gamma_0(23)$ of weight one. It satisfies
$a_2(f)= 1 \in \FF_2$. It can be embedded into weight~$2$ for the same level in two
different ways (multiplying by the Hasse invariant, which does not change the
$q$-expansion, and applying the Frobenius, which sends $q$ to $q^2$).
Consequently, there are two distinct Hecke eigenforms over~$\FF_2$ in weight~$2$
for $\Gamma_0(23)$ whose coefficients at~$2$ are precisely the roots of
$X^2+X+1 \in \FF_2[X]$. The coefficients at the other primes are equal to the
coefficients of~$f$, whence the attached mod~$2$ Galois representations are equal.
Consequently, a converse to Proposition~\ref{prop-use-sturm} cannot exist
(since in this case $m=1$).

\item The trick used in~\cite{CKR} will always work for
deciding whether the representations attached to $f$ and~$g$ are congruent modulo~$\ell^n$:
By applying degeneracy maps at all primes dividing $Nm$ one can force
all coefficients $a_p(f)$ and $a_p(g)$ to be congruent to zero modulo~$\ell^n$ for
all $p \mid Nm$. This allows the application of the Sturm bound. But, usually the
level and hence the bound will be bigger than the bound in
Proposition~\ref{prop-use-sturm}.

\item We mention a point which will be discussed in more detail
in Section~\ref{ex:149times13}. We are mostly interested in congruences of Galois
representations modulo~$\ell^n$ attached to holomorphic eigenforms, hence, it seems
natural to stick to {\em strong} Hecke eigenforms. However, since we formulated
Proposition~\ref{prop-use-sturm} for {\em weak} Hecke eigenforms, we do not need to
have a congruence mod~$\ell^n$ of $\ell$-adic zeros at $p\mid m$, but a simple
equality in the residue ring is enough. Currently, in the algorithm we are not
using this subtle distinction, but, as we will see in the example, it can make a
difference.
\end{enumerate}
\end{rem}

\subsection{Algorithm}\label{sec-alg}

The aim is to study the following problem algorithmically.

\begin{problem}\label{mproblem}
Let $f_1,f_2$ be newforms in levels $N_1,N_2$ and weights $k_1,k_2$.

Determine a finite list of prime powers $\{\ell_1^{n_1},\dots,\ell_r^{n_r}\}$
such that for all $i \in \{1,\dots,r\}$ the $\ell_i$-adic Galois representations
attached to the modular forms $f_1$ and $f_2$ are congruent modulo~$\ell_i^{n_i}$
and are incongruent modulo~$\ell_i^{n_i+1}$,
and for any~$\ell$ distinct from all the~$\ell_i$ the $\ell$-adic Galois
representations of $f_1$ and $f_2$ are incongruent modulo~$\ell$.
\end{problem}

Towards this problem we employ the methods developed in the Section~\ref{sec-cong}.
Due to its greater speed we first apply the congruence number method,
which by Proposition~\ref{ub} gives an upper bound for the possible congruences.
Only if in one of the applications of Corollary~\ref{cor-alg} the upper bound
is unequal to the lower bound we make use of the Newton polygon method.

We hence start by computing the congruence numbers
$c_p=c(P_{f_1,p},P_{f_2,p})$ for all primes $p \nmid N_1 N_2$ up to
some bound (e.g.\ the Sturm bound), where $P_{f_i,p}$ denotes the
characteristic polynomial (in $\ZZ[X]$) of the Hecke operator~$T_p$ acting on the
span of the $\Gal(\Qbar/\QQ)$-conjugacy class~$[f_i]$ of~$f_i$.
Let us number the primes $p_1,p_2,\dots$.
We compute a slightly modified greatest common divisor of all $c_p$,
taking in account only the prime-to-$p$
part of~$c_p$, because we want to disregard the coefficient~$a_p$ when reducing
modulo powers of~$p$.
More precisely, if we have two $c_{p_1}$ and $c_{p_2}$, the first greatest common divisor
that we compute will be
$c=\gcd(c_{p_1}\cdot p_1^{v_{p_1}(c_{p_2})},c_{p_2}\cdot
p_2^{v_{p_2}(c_{p_1})})$.
Once we have one $c$ computed, we can improve it for the next~$p_i$ with
$c'=\gcd(c_{p_i}\cdot p_i^{v_{p_i}(c)},c)$.
The significance of the number~$c'$ is that it gives an upper bound for Problem~$\ref{mproblem}$:
if a prime power $\ell^n$ does not divide~$c'$, then there cannot exist any congruence
modulo~$\ell^n$ between the $\ell$-adic Galois representations attached to $f_1$
and~$f_2$.

Our approach to a solution of Problem~\ref{mproblem} is based on
Theorem~\ref{thm-congr} and Proposition~\ref{prop-use-sturm} in order to
obtain a lower bound, which in favourable cases equals the upper bound~$c'$.
However, whether we use the congruence number method or the Newton
polygon method for computing congruences between zeros of the characteristic
polynomials of the Hecke operators, we have to assume the
following hypothesis, which -- roughly speaking -- says that it is no
loss to work with $P_{f,p}$ instead of with its roots.

\begin{hyp}\label{hyp}
Let $f_1$ and $f_2$ be two newforms and $n \in \NN$.
Suppose that for all primes~$p$ there are embeddings $\sigma_{i,p}: K \hookrightarrow \Qbar$
($i=1,2$) such that
$$\sigma_{1,p}\big(a_p(f_1)\big) \equiv \sigma_{2,p}\big(a_p(f_2)\big) \mod \ell^n.$$
Then there are embeddings $\sigma_1,\sigma_2$ such that
$\sigma_1(f_1) \equiv \sigma_2(f_2) \mod \ell^n$.

An equivalent formulation is the following: If $P_{f_1,p}$ and $P_{f_2,p}$ have roots
congruent modulo~$\ell^n$ (in the sense of Section~\ref{sec-cong}) for all~$p$,
then there are members $\tilde{f}_i$ in the $\Gal(\Qbar/\QQ)$-conjugacy class of~$f_i$
for $i=1,2$ such that $f_1$ is congruent to~$f_2$ modulo~$\ell^n$.
\end{hyp}

In the sequel we shall assume this hypothesis to be satisfied.
Note that by using characteristic polynomials of Hecke operators we lose track
of which form in the $\Gal(\Qbar/\QQ)$-conjugacy class really satisfies a congruence.
By abuse of language we will nevertheless speak of a congruence between
$\rho_{f,\ell}$ and $\rho_{g,\ell}$ modulo~$\ell^n$ when indeed we only
have a congruence of $\rho_{\tilde{f},\ell}$ and $\rho_{\tilde{g},\ell}$ for some
members $\tilde{f}$ and $\tilde{g}$ of the conjugacy classes of $f$ and $g$, respectively.
We now sketch our algorithm for treating Problem~\ref{mproblem}.

\noindent {\bf\underline{Input:}}
$f\in S_k(\Gamma_0(N_f))$ and $g\in S_k(\Gamma_0(N_g))$ be two normalised eigenforms.

\noindent {\bf \underline{Output:}} $(L^-, L^+)$ (for an explanation see below).

\begin{itemize}
\itemsep=0cm plus 0pt minus 0pt
\item (Upper bound) For every prime $p\nmid  N_f N_g$ up to the Sturm bound $B$ (see Theorem~\ref{thm-fin-gen}),
we compute the congruence number $c_p=c(P_{f,p},P_{g,p})$ and we
calculate $L^+=\gcd_{p\le B}(c_p)$ with the
modified greatest common divisor described above.
We recall that $P_{f,p}$ denotes the characteristic
polynomial of the Hecke operator~$T_p$ acting on the span~$[f]$ of the Galois
conjugacy class of~$f$, which can for instance be obtained as the characteristic polynomial
of the action of~$T_p$ on a suitable modular symbols space.

\item For every $\ell\mid L^+$, we compute $L^-_{1,\ell}=\min_{p\le B}(\ell^{d_p})$,
where $\ell^{d_p}$ is the maximal power of~$\ell$ modulo which
$P_{f,p}$ and $P_{g,p}$ have a root in common. This number is
obtained from the congruence number method
if the value returned by it is best possible, i.e.\ if we are in case
(\ref{part-nice}) or (\ref{part-one}) of Corollary~\ref{cor-alg}.
Otherwise, the Newton polygon method is employed.
We then form the product
$L^-_1=\prod_{\ell\mid L^+}L^-_{1,\ell}$.

\item Suppose for this step that $N_g = mN_f$
and that $\rhobar_{f,\ell}$ and $\rhobar_{g,\ell}$ are absolutely irreducible.
Then, for every $\ell\mid L^+$
such that $v_\ell(L^+)\ne v_\ell(L^-_1)$,
we compute $L^-_{2,\ell}=\min_{p\le B}(\ell^{\tilde d_p})$
as follows:
If $p \nmid m$, then we put $\tilde{d}_p = d_p$.
If $p \mid m$, we let $\ell^{\tilde{d}_p}$ be the maximal power of~$\ell$ modulo which
$\tilde{P}_{f,p}$ and $P_{g,p}$ have a root in common
with $\tilde{P}_{f,p}$ as in Equation~(\ref{charpoly}). This number is again
calculated by the congruence number method or the Newton polygon method as in
the previous step.
Again we compute $L^-_2=\prod_{\ell\mid L^+}L^-_{2,\ell}$.

\item We compute $L^-=\prod_{\ell\mid L^+}\max(L^-_{1,\ell},L^-_{2,\ell})$.

\item Return $(L^-,L^+)$.
\end{itemize}

Proposition~\ref{ub} ensures that $L^+$ is an upper bound, i.e.\ that $\rho_{f,\ell}$ and
$\rho_{g,\ell} $ are incongruent modulo~$\ell^m$
(more precisely, this holds for any members of the conjugacy classes of~$f$ and~$g$)
if $\ell^m \nmid L^+$.
Theorem~\ref{thm-congr} guarantees that $L^-_1$ is a lower bound
(under Hypothesis~\ref{hyp}), meaning that under the hypothesis
$\rho_{f,\ell}$ and $\rho_{g,\ell} $ are congruent modulo~$\ell^n$ if $\ell^n \mid L_1^-$
(with the slight abuse of language pointed out above).
The lower bound $L_1^-$ will in general be very bad (e.g.\ $1$) due to the Hecke
operators $T_p$ for $p \mid m$ (in the situation of the third step). This is taken
care of in the third step and Proposition~\ref{prop-use-sturm} tells us that
$L_2^-$ is a lower bound in the same sense as before (still under Hypothesis~\ref{hyp}).
Consequently, $L^-$ is a lower bound under Hypothesis~\ref{hyp}.

\begin{rem}
We point out that this algorithm might miss a congruence modulo~$\ell^n$ due to
the Hecke operator~$T_\ell$. Hence, one might want to exclude the operators $T_\ell$ in
all the steps. Then, however, we do not have the congruence of $g$ with an oldform of~$f$
(as in Proposition~\ref{prop-use-sturm}), hence, the congruence of the Galois representations
suggested by the output of the algorithm will not be a proved result even
under Hypothesis~\ref{hyp} (but the correct one in most cases).
\end{rem}

\section{Examples and numerical data}\label{sec-examples}

In this section we present some cases which were computed using
the algorithm described above and which we consider interesting.
Several more examples can be found in~\cite{T}.
For our calculations we used the computer algebra system {\sc Magma} (\cite{Magma}).

\subsection{Examples of congruences in the same level}

\begin{table}[t]
\begin{center}
\begin{tabular}{|l|l||l|l||l|l|}
\hline
$N_1$ & $i_1$ & $N_2$ & $i_2$ & lower bound & upper bound \\
\hline

$1479$ & $16 $ & $1479$&$ 8$ & $ 2^7$ & $ 2^7$\\
$1027$ & $2$ & $1027$ & $1$ & $2^5$ & $2^5$\\
$602$ & $8$ & $602$ & $7$ & $2^5$ & $2^5$\\
$1454$ & $7$ & $1454$ & $1$ & $3^4$ & $3^4$\\
$1171$ & $4$ & $1171$ & $2$ & $3^4$ & $3^4$\\
$1147$ & $6$ & $1147$ & $5$ & $7^3$ & $7^3$\\
$1726$ & $6$ & $1726$ & $3$ & $5^3$ & $5^3$\\
$1629$ & $4$ & $1629$ & $3$ & $5^3$ & $5^3$\\
$613$ & $2$ & $613$ & $1$ & $7\cdot47^2$ & $7\cdot47^2$\\
$1939$ & $4$ & $1939$ & $2$ & $37^2\cdot4423$ & $37^2\cdot4423$\\
$1906$ & $5$ & $1906$ & $3$ & $19^2$ & $19^2$\\
$1763$ & $8$ & $1763$ & $5$ & $3\cdot13^2$ & $3\cdot13^2$\\
$1761$ & $8$ & $1761$ & $7$ & $2\cdot8581981$ & $2\cdot8581981$\\
$1241$ & $2$ & $1241$ & $1$ & $1933\cdot8713$ & $1933\cdot8713$\\

$ 71 $ & $ 2 $ & $71$ & $ 1$ & $ 2\cdot 3^2 $&$2 \cdot 3^2$\\
$ 109$ &$ 3 $&$109 $&$1$&$2^2$&$ 2^2 $\\
$ 155$&$ 4 $&$155$&$ 2$&$ 2^4 $&$2^4 $\\
$ 233 $&$3 $&$233 $&$1$&$ 3^3$&$ 3^3 $\\
$ 785 $&$2 $&$785 $&$1 $&$7^3 $&$7^3 $\\
$ 1073 $&$6 $&$1073 $&$3 $&$2\cdot 17^2  $&$ 2\cdot 17^2$\\
$ 1481 $&$3 $&$1481 $&$1 $&$5^2 \cdot 2833 $&$ 5^2 \cdot 2833 $\\
\hline
\end{tabular}
\caption{Extract from the computational results.}
\label{table-one}
\end{center}
\end{table}

We computed all congruences between modular forms of weight~$2$
and the same level up to level $2000$.
In Table~\ref{table-one},
$(N_j,i_j)$ means the $i_j$-th form in level~$N_j$ for $j=1,2$
(according to an internal ordering in {\sc Magma}), where
in these cases we have $N_1=N_2$.
In all these cases, we found $L^- = L^+$ so that under Hypothesis~\ref{hyp}
we obtained all congruences.

\begin{itemize}
\item The biggest exponents that we found appear in $2^7$ and $2^5$.
\item For $n=4$, we find some congruences modulo $3^4$ (also
modulo $2^4$).
\item For $n=3$, the primes $\ell=5$ and $\ell=7$ appear.
\item For $n=2$ we already have many different primes, $47^2$ being
the biggest square of a prime that we found.
\item For $n=1$ we just listed some of the biggest congruences that we
found. $2\cdot8581981=17163962$ and $1933\cdot8713 = 16842229$ are
just two examples of congruences, but in this case we had several
primes to choose from.
\end{itemize}

\subsection{Simple example for strong $\neq$ weak}\label{sec-weakstrong}

We now analyse the example with the smallest level in the above table
more thoroughly.
On $\Gamma_0(71)$ there are two $\Gal(\Qbar/\QQ)$-conjugacy classes of newforms
in weight~$2$. The coefficient fields of both of them are isomorphic; they have
degree~$3$, discriminant~$257$ and are non-Galois. The prime~$3$ factors in two
prime ideals~$\fP_1$ and $\fP_2$ of residue degrees~$1$ and~$2$. This means that
each of the two $\Gal(\Qbar/\QQ)$-conjugacy classes gives us precisely one
strong Hecke eigenform~$f_i$ modulo~$3^n$ with coefficients in~$\ZZ/3^n\ZZ$ for $i=1,2$;
the others taken modulo~$3$ have coefficients in~$\FF_9$.

We compute that $f_1$ and $f_2$ are congruent modulo~$9$, but incongruent modulo~$27$.
Let $\TT \subset \End_\CC(S_2(\Gamma_0(71)))$ be the Hecke algebra, i.e.\ the
subring generated by the Hecke operators. The above discussion shows that there
is a maximal ideal~$\fm$ of $\hat{\TT} := \TT \otimes_\ZZ \ZZ_3$ such that
the localisation $\hat{\TT}_\fm$ has two minimal prime ideals, corresponding to
the two strong Hecke eigenforms $f_1$ and~$f_2$. A computer calculation yields
that $\hat{\TT}_\fm \otimes_{\ZZ_3} \ZZ/9\ZZ \cong \ZZ/9\ZZ [X]/(X^2)$.
Thus, we have three weak Hecke eigenforms modulo~$9$ coming from $\hat{\TT}_\fm$, namely
$$ \hat{\TT}_\fm \twoheadrightarrow \hat{\TT}_\fm \otimes_{\ZZ_3} \ZZ/9\ZZ
\cong \ZZ/9\ZZ [X]/(X^2)
\xrightarrow{X \mapsto 0 \textrm{ or } X \mapsto 3 \textrm{ or } X \mapsto 6} \ZZ/9\ZZ.$$
Since we know that there is only one strong Hecke eigenform modulo~$9$, two of them
cannot be strong.

\subsection{Example in levels $149$ and $149 \cdot 13$}\label{ex:149times13}

On $\Gamma_0(149)$ for weight~$2$ there are two $\Gal(\Qbar/\QQ)$-conjugacy
classes of newforms. The degrees of the coefficient fields are $3$ and $9$.
Let $f$ be any of the forms whose coefficient field $\QQ_f$ has degree~$9$.
The prime $3$ is unramified in~$\QQ_f$ and
there is a prime $\fP$ of residue degree~$1$ in the ring of integers $\cO_f$ of~$\QQ_f$.

Mazur's Eisenstein ideal (\cite{Eisenstein})
shows that the residual representation $\rhobar_{f,\fP}$
of~$f$ modulo~$\fP$ is irreducible, since $149$ is a prime number and
$3$ does not divide $149-1$.
We first want to determine the image of the residual representation.
A quick computation of a couple of coefficients of~$f$ shows that
the image of $\rhobar_{f,\fP}$ contains all possible combinations of trace and determinant.
Consulting the list of subgroups of $\GL_2(\FF_3)$ tells us that next to
the full $\GL_2(\FF_3)$ there is only one other subgroup satisfying this property.
That subgroup, however, does not contain any element of order~$3$. Due to the
semistability at $13$ and~$149$ this group is excluded, whence
the image is the full $\GL_2(\FF_3)$.

There is a newform~$g$ of weight~$2$ on $\Gamma_0(13\cdot 149)$ and a prime ideal
$\Lambda$ dividing~$3$ in its coefficient field such that the strong Hecke eigenform
of~$g$ obtained by reducing its $q$-expansion modulo~$\Lambda$ is equal to
the strong Hecke eigenform of $f$ modulo~$\fP$ at all prime coefficients except at~$13$.
In fact, our algorithm gives us a congruence modulo~$3^{10}$ (in the sense defined before)
at all primes up to the Sturm bound, except~$13$. Moreover, $3^{10}$ is also an upper bound.
At the prime~$13$ we want to apply Proposition~\ref{prop-use-sturm} (i.e.\ the third
item of the algorithm), and we hence apply the methods from Corollary~\ref{cor-alg}
to $P_{g,13}$ and $\tilde{P}_{f,13}$. However, the upper and the lower bounds we
obtain with this method are $3^9$. Hence, the output of our algorithm would be a
congruence modulo~$3^9$ of the Galois representations attached to $f$ and~$g$
as lower bound and $3^{10}$ as upper bound.
We analyse the situation a bit more closely by hand. The polynomial $P_{g,13}$
is equal to $(X+1)^{80}$. The polynomial $\tilde{P}_{f,13} = Q^2$ with
$Q \in \ZZ[X]$ an irreducible polynomial of degree~$18$.
Evaluating $Q$ at $-1$ (the zero of~$P_{g,13}$) gives
$2^6 \cdot 3^{10} \cdot 6869$. This means that there is a {\em weak} Hecke
eigenform $\tilde{f}$ in the oldspace of~$f$ modulo~$3^{10}$ such that
$\tilde{f}(T_{13}) = -1$.
Hence, Proposition~\ref{prop-use-sturm} yields that $\tilde{f}$ and $g$
are congruent modulo~$3^{10}$ as weak Hecke eigenforms. Consequently, the
attached Galois representations of $f$ and~$g$ are congruent modulo~$3^{10}$.

We give a more formal argument for the existence of the weak Hecke eigenform
modulo~$3^{10}$. Let $\TT$ be the Hecke algebra
on $S_2(\Gamma_0(149\cdot 13))$ (as $\ZZ$-algebra) and let $\TT_{[f]}^\old$ be the
Hecke algebra (as $\ZZ$-algebra) on the image of~$[f]$ under the $13$-degeneracy map,
where as before $[f]$ denotes the span of the Galois conjugacy classes of~$f$. By restricting
Hecke operators, we obtain a surjective ring homomorphism
$\TT \twoheadrightarrow \TT_{[f]}^\old$. The algebra $\TT_{[f]}^\old$ is generated
by the identity matrix and~$\tilde{T}_{13}$ (see Equation~(\ref{bigmatrix})).
Since the minimal polynomial of $\tilde{T}_{13}$
is either $Q$ or $Q^2$, the composition
$$ \TT \twoheadrightarrow \TT_{[f]}^\old \xrightarrow{\tilde{T}_{13} \mapsto -1} \ZZ/3^{10}\ZZ$$
is a well-defined ring homomorphism, i.e.\ the desired weak Hecke eigenform modulo~$3^{10}$.

\subsection{Congruences with Eisenstein series modulo~$\ell^n$}

Let $f \in S_2(\Gamma_0(N))$ such that $\rhobar_{f,\ell}$ is reducible
(and semi-simple by definition).
This means that $f$ is congruent modulo~$\ell$ to an Eisenstein series
in the same level and weight at almost all primes. The converse of this statement
also holds.
In the context of this article,
it is natural to study congruences between newforms and Eisenstein series modulo~$\ell^n$
and to do so via the congruence number and the Newton polygon method.
By computing congruences modulo~$\ell^n$ with Eisenstein series,
we study up to which $\ell^n$ the representation $\rhobar_{f,\ell^n}$
has the same traces at the first couple of Frobenius elements at good primes
as an extension of the cyclotomic character modulo~$\ell^n$ by the trivial representation.

Let $f$ be a newform of weight~$k$ and level~$N$.
We implemented an algorithm, which for all primes~$p \nmid N$ up to the Sturm bound
computes the maximal prime powers modulo which $P_{f,p}$
(as before, this is the characteristic polynomial of $T_p$ acting on~$[f]$)
and the characteristic polynomial of~$T_p$ acting on the Eisenstein subspace in
the given level and weight have a root in common.
We then proceed as earlier, obtaining an upper bound
for a congruence with an Eisenstein series as well as an
unproved lower bound (note that we do not take all operators into account).

A famous theorem of Mazur's (\cite{Eisenstein}) states that in weight~$2$
and prime level~$N$ there is a cusp form
which is congruent to the Eisenstein series modulo~$\ell$ at almost all primes
for every~$\ell$ dividing the numerator of $\frac{N-1}{12}$.
One can ask in how far this theorem holds modulo~$\ell^n$. It quickly turns out
that a too naive generalisation is false.
We propose to study the following in a subsequent paper.
Let $f_1,\dots,f_r$ be all newforms in prime level~$N$ and weight~$2$ for the trivial
Dirichlet character. For $i=1,\dots,r$ let $\ell^{n_i}$ be the highest power of~$\ell$
such that $f_i$ is congruent at almost all primes to the Eisenstein series of level~$N$ and weight~$2$
modulo~$\ell^{n_i}$. Put $n := n_1 + \ldots + n_r$.

\begin{question}
Is $n$ at least as big as (or even equal to) the $\ell$-valuation of the numerator of $\frac{N-1}{12}$?
\end{question}

\subsection{Level raising modulo~$\ell^n$}

Let $f \in S_2(\Gamma_0(N))$ be a newform. The term {\em level
raising modulo~$\ell^n$} in the simplest case refers to the problem
of identifying primes~$p\nmid N$ such that there is a newform~$g$ in
$S_2(\Gamma_0(Np))$ with the property that $f$ and $g$ are congruent
modulo~$\ell^n$ at almost all primes. A necessary condition for
level raising of the form $f$ modulo~$\ell$ at the prime~$p \nmid N$
when its Galois representation is residually irreducible, is that $\ell$
divides the congruence number $c(P_{f,p},X-(p+1))$ or the congruence
number $c(P_{f,p},X+(p+1))$. It is a famous theorem of Ribet's
(\cite{R}) that the converse also holds (modulo~$\ell$).

It is natural to ask whether or in which sense level raising
generalises to congruences modulo~$\ell^n$. We start by an
observation which we consider very interesting.
Let $f$ be the only newform on $\Gamma_0(17)$ in weight~$2$ and let $p=59$.
The coefficient
$a_{59}(f) = -12$ and we find that $9$ divides $c(P_{f,59},X-60) =
c(X+12,X-60)=72$ and that $3$ divides $c(P_{f,59},X+60) =
c(X+12,X+60)=48$. However, there does not seem to be a congruence
modulo~$9$ of $f$ with any form in level $17 \cdot 59$. Instead,
there appear to be three newforms in that level which are congruent
to~$f$ modulo~$3$ at almost all primes. Hence, we conclude that the
condition that $\ell^n$ divides one of the above congruence numbers
is not a sufficient one for level raising of strong Hecke
eigenforms. This confirms a remark by Richard Taylor.\footnote{This
remark was made in the Problem Book for the MSRI Modular Forms
Summer Workshop organised by William Stein in 2006.}

We propose to study the following question in a subsequent paper.
Let $f \in S_2(\Gamma_0(N))$ be some newform and let $p \nmid N$ be a prime.
Further, let $g_1,\dots,g_r$ be all newforms in $S_2(\Gamma_0(Np))$.
For $i=1,\dots,r$ let $\ell^{n_i}$ be the highest power of~$\ell$
such that $g_i$ is congruent to~$f$ modulo~$\ell^{n_i}$ at almost all primes.
Put $n := n_1 + \ldots + n_r$ and let $c$ be the maximum integer
such that $P_{f,p}$ and $X^2-(p+1)^2$ have a root in common modulo~$\ell^c$.

\begin{question}
Is $n$ equal to the $\ell$-valuation of~$c$?
\end{question}

An inequality (in a greater generality) is
provided by Theorem~2 of~\cite{D}.

\end{document}